\documentclass[a4paper,12pt]{article}

\usepackage{graphicx}
\usepackage{mathptmx}

\usepackage{amssymb,amsfonts,amsmath,epsfig,float,graphicx}

\newcommand{\be}{\begin{eqnarray}}
\newcommand{\ee}{\end{eqnarray}}
\newcommand{\nnee}{\nonumber\ee}

\newcommand{\upd}{{\rm \,d}}

\newcommand{\beginproof}{\par\strut\vskip 0.1cm\noindent{\bf Proof}\par}
\newcommand{\beginproofof}[1]{\par\strut\vskip 0.1cm\noindent{\bf Proof of #1}\par}

\def\qed{\par\strut\hfill$\square$\par\vskip 0.5cm}

\newtheorem{theorem}{Theorem}[section]
\newtheorem{proposition}[theorem]{Proposition}
\newtheorem{lemma}[theorem]{Lemma}

\newcommand{\tinyc}{{\mbox{\tiny c}}}

\newcommand{\Ro}{{\mathbb R}}
\newcommand{\Eo}{{\mathbb E}}
\newcommand{\Io}{{\mathbb I}}

\newcommand{\emp}{{\mbox{\tiny emp}}}
\newcommand{\pemp}{p^\emp}
\newcommand{\Eoemp}{\Eo^\emp}

\newcommand{\Ac}{A^{\mbox{\tiny c}}}
\newcommand{\Bc}{B^{\mbox{\tiny c}}}

\title{Update of prior probabilities by minimal divergence} 

\author{Jan Naudts\\ Universiteit Antwerpen\\
 Antwerpen, Belgium\\
}

\begin{document}
\maketitle

\begin{abstract}
The present paper investigates
the update of an empirical probability distribution with the results of a new set of observations.
The optimal update is obtained by minimizing either the Hellinger distance or
the quadratic Bregman divergence. The results obtained by the two methods differ.
Updates with information about conditional probabilities are considered as well.
%

\end{abstract}

\section{Introduction}
\label{intro}

In the present work prior probabilities are assumed to be known.
The approach is then to look for updated probabilities
that reproduce the observed probabilities
while taking the prior into account.

No further model assumptions are imposed. 
Hence, the statistical model under consideration consists of all probability distributions
that are consistent with the newly obtained empirical data.
Internal consistency of the empirical data ensures that the model is not empty.
The update is the model point that minimizes the chosen divergence function
from the prior to the manifold of the model.

In the context of Maximum Likelihood Estimation (MLE) the model is usually well-known,
the dimension of the model can be kept low 
and properties of the model can be used to ease the calculations.
It can then happen that the model is misspecified \cite{WH82} and that the update
is only a good approximation of the empirical data. Here the model is dictated
by the newly acquired empirical data and the update is forced to reproduce the measured data.

In Bayesian statistics the update $q(B)$ of the probability $p(B)$ of an event $B$ equals
\be
q(B)&=&p^\emp(A)\,p(B|A)+p^\emp(A^\tinyc)\,p(B|A^\tinyc).
\nnee
The quantities $p^\emp(A)$ and $p^\emp(A^\tinyc)$ are the empirical probabilities obtained
after repeated measurement of event $A$ and its complement $A^\tinyc$ .
This result is obtained also when minimizing the Hellinger distance between the prior 
and the model manifold. A proof of the latter follows later on in Section \ref {sect:optimal}.

One incentive for starting the present work is a paper of
Banerjee, Guo, and Wang \cite{BGW05,FSG08}.
They consider the problem of predicting a random variable $Z_1$
given observations of a random variable $Z_2$.
It is well-known that the conditional expectation, as defined by Kolmogorov, 
is the optimal predictor.
They show that this statement remains true when the metric distance
is replaced by a Bregman divergence. 
Their result is limited to the updated
expectation of random variables $Z_1$ that are
functions of the measured random variable $Z_2$.
It is shown in Theorem \ref {theorem2}
below that a proof in a more general context yields a deviating result.

The present work is inspired by the current practices in Information Geometry
\cite{AN00,AS16,AJLS18} where minimization of divergences is an important tool.
In Statistical Physics a divergence is called a relative entropy.
Its importance was noted rather late in the twentieth century, after the work of
Jaynes on the maximal entropy principle \cite{JE57}. 
Estimation in the presence of hidden variables by minimizing a divergence function
is briefly discussed in Chapter 8 of \cite{AS16}.

The next Section fixes notations.
Section \ref{sect:geom} collects some results about the squared Hellinger distance
and the quadratic Bregman divergence.
Section \ref {sect:optimal} discusses the optimal choice and 
contains the Theorems \ref{theorem1} and \ref {theorem2}.
The proof of the theorems can be adapted to cover the situation that
a subsequent measurement also yields information on conditional probabilities.
This is done in Section \ref {sect:opt:hell}.
A final section summarizes the results of the paper.

\section{Empirical data}
\label {sect:emp}

Consider a probability space $X,\mu$.
A measurable subset $A$ of $X$ is called an {\sl event}.
Its probability is denoted $p(A)$ and is given by
\be
p(A)&=&\int_X\Io_A(x)\upd\mu(x),
\nnee
where $\Io_A(x)$ equals 1 when $x\in A$ and 0 otherwise. 
The conditional expectation of a random variable $f$ given an event $A$
with non-vanishing probability $p(A)$ is given by
\be
\Eo_\mu f|A&=&\frac 1{p(A)}\Eo_\mu f\,\Io_A.
\nnee

The probability space $X,\mu$ reflects the prior knowledge of the system at hand.
When new data become available an update procedure is used to select the posterior
probability space. The latter is denoted $X,\nu$ in what follows.
The corresponding probability of an
event $A$ is denoted $q(A)$.

The outcome of repeated experiments is the empirical probability distribution of the events,
denoted $\pemp(A)$.
The question at hand is then to establish a criterion for finding
the update $\nu$ of the probability distribution $\mu$ that is as close as possible to
$\mu$ while reproducing the empirical results.

The event $A$ defines a partition $A,\Ac$ of the probability space $X,\mu$.
As before $\Ac$ denotes the complement of $A$ in $X$.
In what follows a slightly more general situation is considered
in which the event $A$ is replaced by a partition $(O_i)_{i=1}^n$ of the
measure space $X,\mu$ into subsets with non-vanishing probability.
The notations $p_i$ and $\mu_i$ are used, with
\be
p_i=p(O_i)
\quad\mbox{ and }\quad
\upd\mu_i(x)=\frac{1}{p_i}\Io_{O_i}(x)\upd \mu(x).
\label{update:muidef}
\ee

Introduce the random variable $g$ defined by $g(x)=i$ when $x\in O_i$.
Repeated measurement of the random variable $g$ yields the empirical probabilities
\be
\pemp_i&=&\mbox{Emp Prob }\{x:\,g(x)=i\}.
\nnee
They may deviate from the prior probabilities $p_i$.
In some cases one also measures the conditional probabilities
\be
\pemp(B|O_i)&=&\mbox{Emp Prob of }B\mbox{ given that }g(x)=i
\nnee
of some other event $B$.

\section{A geometric approach}
\label{sect:geom}

In this section two divergences are reviewed, the squared Hellinger distance
and the quadratic Bregman divergence.

\subsection{Squared Hellinger distance}

For simplicity the present section is restricted to the case that the sample space $X$ is the real line.

Given two probability measures $\mu$ and $\sigma$, both absolutely continuous
w.r.t.~the Lebesgue measure,
the squared Hellinger distance is the divergence $D_2(\sigma||\mu)$ defined by
\be
D_2(\sigma||\mu)
&=&
\frac 12\int_\Ro \left(\sqrt{\frac{\upd \sigma}{\upd x}}-\sqrt{\frac{\upd \mu}{\upd x}}\right)^2\upd x.
\nnee
It satisfies
\be
D_2(\sigma||\mu)
&=&
1-\int_\Ro\sqrt{\frac{\upd \sigma}{\upd x}\frac{\upd \mu}{\upd x}}\upd x.
\nnee

Let $(O_i)_i$ be a partition of $X,\mu$
and let $g(x)=i$ when $x$ belongs to $O_i$, as before.
Let $p_i$ and $\mu_i$ be defined by (\ref {update:muidef}).
Consider the following functions of $i$, with $i$ in $\{1,\cdots,n\}$,
\be
\tau^{(1)}(i)&=&\mu,\quad\mbox{ independent of }i,\cr
\tau^{(2)}(i)&=&\mu_i,\cr
\tau^{(3)}(i)&=&\sigma_i,
\nnee
where each of the $\sigma_i$ is a probability distribution with support in $O_i$.
The empirical expectation of a function $f(i)$ is given by
$\Eoemp f=\sum_i\pemp_i f(i)$.

\begin{proposition}
\label {prop:minemp}
If $\pemp_i>0$ for all $i$ and $\sum_i\pemp_i=1$ then one has
\be
\Eoemp D_2(\tau^{(1)}|| \tau^{(3)})\ge \Eoemp D_2(\tau^{(1)}|| \tau^{(2)})
\nnee
with equality if and only if $\sigma_i=\mu_i$ for all $i$.
\end{proposition}

First prove the following two lemmas.

\begin{lemma}
Assume that the probability measure $\nu_i$
is absolutely continuous w.r.t. the measure $\mu_i$, with Radon-Nikodym derivative
given by $\upd\nu_i(x)=f_i(x)\upd\mu_i$.
Then one has
\be
D_2(\mu|| \sigma_i)-D_2(\mu||\nu_i)
&=&
\sqrt{p_i}\left[D_2(\mu_i|| \sigma_i)-D_2(\mu_i||\nu_i)\right]
\nnee
and
\be
D_2(\mu_i||\nu_i)&=&1-\int_{O_i}\sqrt{f_i(x)}\upd\mu_i(x).
\nnee

\end{lemma}

\beginproof

One calculates
\be
D_2(\mu|| \sigma_i)-D_2(\mu||\nu_i)
&=&
\int_{\Ro}\sqrt{\frac{\upd\mu}{\upd x}}\left[
\sqrt{\frac{\upd\nu_i}{\upd x}}-\sqrt{\frac{\upd\sigma_i}{\upd x}}
\right]\upd x\cr
&=&
\sqrt{p_i}
\int_{O_i}\sqrt{\frac{\upd\mu_i}{\upd x}}\left[
\sqrt{\frac{\upd\nu_i}{\upd x}}-\sqrt{\frac{\upd\sigma_i}{\upd x}}
\right]\upd x\cr
&=&
\sqrt{p_i}\left[\int_{O_i}\sqrt{f_i(x)}\upd\mu_i(x)
-\int_{O_i}\left[\frac{\upd\mu_i}{\upd x}\frac{\upd\sigma_i}{\upd x}\right]^{1/2}\upd x\right]\cr
&=&
\sqrt{p_i}\left[
\int_{O_i}\sqrt{f_i(x)}\upd\mu_i(x)-1+
D_2(\mu_i|| \sigma_i)\right].
\nnee
Now take $\sigma_i=\nu_i$ to obtain the desired results.
\qed

\begin{lemma}
\label{prop:pythhel}
(Pythagorean relation) For any $i$ is
\be
D_2(\mu|| \sigma_i)
&=&
D_2(\mu|| \mu_i)+\sqrt{p_i}D_2(\mu_i||\sigma_i).
\nnee
\end{lemma}

\beginproof
The proof follows by taking $\nu_i=\mu_i$ in the previous lemma.
\qed

\beginproofof{Proposition \ref  {prop:minemp}.}

From the previous lemma it follows that 
$D_2(\tau^{(1)}|| \tau^{(3)})\ge D_2(\tau^{(1)}|| \tau^{(2)})$.
Note that $\sigma_i=\mu_i$ implies that $\tau^{(3)}=\tau^{(2)}$
and hence $D_2(\tau^{(1)}|| \tau^{(3)})= D_2(\tau^{(1)}|| \tau^{(2)})$.
Conversely, if 
\be
\Eoemp D_2(\tau^{(1)}|| \tau^{(3)})&=& \Eoemp D_2(\tau^{(1)}|| \tau^{(2)})
\nnee
then it follows from the previous proposition that $\Eoemp D_2(\tau^{(2)}|| \tau^{(3)})=0$.
If in addition $\pemp_i>0$ for all $i$ then it follows that for all $i$
\be
0&=&D_2(\tau^{(2)}(i)|| \tau^{(3)}(i)).
\nnee
Because the squared Hellinger distance is a divergence this implies that
$\tau^{(2)}(i)=\tau^{(3)}(i)$, which is equivalent with $\mu_i=\sigma_i$.
\qed

\subsection{Bregman divergence}

In the present section the squared Hellinger distance, which is an f-divergence,
is replaced by a divergence of the Bregman type. In addition let $X$ be a finite set.
Then there exists for each of the elements $O_i$ of the partition of $X$ a counting measure
$\rho_i$ such that
\be
\rho_i(x)&=&\frac{1}{|O_i|}
\quad\mbox{ if } x\in O_i,\cr
&=&0
\quad\mbox{ otherwise}.
\label{breg:countdef}
\ee

Fix a strictly convex function $\phi:\,\Ro\mapsto\Ro$.
The Bregman divergence of the probability measures $\sigma$ and $\mu$
is defined by
\be
D_\phi(\sigma||\mu)
&=&
\sum_x\left[\phi\left(\sigma(x)\right)-\phi\left(\mu(x)\right)
-\left(\sigma(x)-\mu(x)\right)
\phi'\left(\mu(x)\right)\right]
\nnee
In the case that $\phi(x)=x^2/2$, which is used below, it becomes
\be
D_\phi(\sigma||\mu)
&=&
\frac 12\sum_x\left[\sigma(x)-\mu(x)\right]^2.
\label{breg:2def}
\ee
For convenience, this case is referred to as the {\em quadratic Bregman divergence}.

The following result, obtained with the quadratic Bregman divergence, is 
more elegant than the result of Lemma \ref {prop:pythhel}.

\begin{proposition}
\label{prop:bregpyth}
Consider the quadratic Bregman divergence $D_\phi$ as given by (\ref {breg:2def}).
Let $\nu_i=p_i\mu_i+(1-p_i)\rho_i$. 
Let $\sigma_i$ be any probability measure with support in $O_i$.
Then the following Pythagorean relation holds.
\be
D_\phi(\mu|| \sigma_i)
&=&
D_\phi(\mu||\nu_i)+D_\phi(\nu_i||\sigma_i).
\nnee

\end{proposition}

\beginproof
One calculates
\be
D_\phi(\mu|| \sigma_i)-D_\phi(\mu||\nu_i)
&=&
D_\phi(\nu_i||\sigma_i)
+\sum_x\left[\mu(x)-\nu_i(x)\right]\,
\left[\phi'\left(\nu_i(x)\right)
-\phi'\left(\sigma_i(x)\right)\right]\cr
&=&
D_\phi(\nu_i||\sigma_i)
+\sum_{x\in O_i}\left[p_i\mu_i(x)-\nu_i(x)\right]\,
\left[\phi'\left(\nu_i(x)\right)
-\phi'\left(\sigma_i(x)\right)\right]\cr
&=&
D_\phi(\nu_i||\sigma_i)
-(1-p_i)\frac{1}{|O_i|}\sum_{x\in O_i}\left[\phi'\left(\nu_i(x)\right)
-\phi'\left(\sigma_i(x)\right)\right].
\nnee
Use now that $\phi'(u)=u$ and the normalization of the probability measures $\nu_i$
and $\sigma_i$ to find the desired result.

\qed

\section{The optimal choice}
\label{sect:optimal}

\subsection{Updated probabilities}

The following result proves that the standard Kolmogorovian definition of the conditional probability
minimizes the Hellinger distance between the prior probability measure $\mu$ and the 
updated probability measure $\nu$. The optimal choice of the updated probability measure $\nu$
is given by corresponding probabilities $q(B)$. They satisfy
\be
q(B)&=&\sum_{i=1}^n\pemp_i p(B|O_i)
\quad\mbox{ for any event }B.
\nnee

\begin{theorem}
 \label{theorem1}
Let be given a partition $(O_i)_{i=1}^n$ of the probability space $X,\mu$ with $X=\Ro$.
Let $\mu_i$ be given by (\ref {update:muidef}).
Let $p_i=p(O_i)>0$ denote the probability of the event $O_i$
and let be given strictly positive empirical probabilities $\pemp_i$, $i=1,\cdots, n$. 
The squared Hellinger distance $D_2(\sigma||\mu)$ as a function of $\sigma$
is minimal if and only if $\sigma_i=\mu_i$ for all $i$.
Here, $\sigma$ is any probability measure on $X$ satisfying 
\be
\sigma=\sum_{i=1}^n\pemp_i\sigma_i,
\nnee
and each of the $\sigma_i$ is a probability measure with support in $O_i$
and absolutely continuous w.r.t.~$\mu_i$.

\end{theorem}

Note that the probability measure $\nu$ given by
\be
\nu(x)&=&\sum_{i=1}^n\pemp_i\mu_i(x)
\nnee
uses the Kolmogorovian conditional probability as the predictor because 
the probabilities determined by the $\mu_i$ are obtained
from the prior probability distribution $\mu$ by $p_i(x)=p(x|O_i)$.
By the above theorem this predictor is the optimal one w.r.t.~the squared Hellinger distance.

\beginproof
With the notations of the previous section is
\be
D_2(\sigma||\mu)&=&\Eoemp D_2(\tau^{(1)}|| \tau^{(3)}).
\nnee
Proposition \ref{prop:minemp} shows that it is minimal if and only if $\sigma_i=\mu_i$
for all $i$.

\qed

\def\Dav{D^{\mbox{\tiny av}}}

Next, consider the use of the quadratic Bregman divergence
in the context of a finite probability space.

\begin{theorem}
 \label{theorem2}
Let be given a partition $(O_i)_{i=1}^n$ of the finite probability space $X,\mu$.
Let $\rho_i$ be the counting measure on $O_i$ defined by (\ref {breg:countdef}).
Let $\mu_i$ be given by (\ref {update:muidef}).
Let $p_i=p(O_i)>0$ denote the probability of the event $O_i$
and let be given strictly positive empirical probabilities $\pemp_i$, $i=1,\cdots, n$
summing up to 1. 
Assume that
\be
\pemp_i\ge p_i\left[1-|O_i|\mu_i(x)\right]
\quad\mbox{ for all }x\in O_i\mbox{ and for }i=1,\cdots,n.
\label{opt:theo2:assum}
\ee
Then the following hold.
\begin{description}
\item [1)\,] A probability distribution $\nu$ is defined by $\nu=\sum_i\pemp_i\nu_i$ with
\be
\nu_i&=&\left(1-\frac{p_i}{\pemp_i}\right)\rho_i+\frac{p_i}{\pemp_i}\mu_i.
\nnee
\item [2)\,] Let $\sigma$ be any probability measure on $X$ satisfying 
$\sigma=\sum_{i=1}^n\pemp_i\sigma_i$, where each of the $\sigma_i$
is a probability distribution with support in $O_i$. Then the quadratic Bregman divergence 
satisfies the Pythagorean relation
\be
D_\phi(\sigma||\mu)&=&D_\phi(\nu||\mu)+\sum_{i=1}^n(\pemp_i)^2D_\phi(\sigma_i||\nu_i).
\label{opt:pythbreg}
\ee
\item [3)\,] The quadratic Bregman divergence $D_\phi(\sigma||\mu)$
is minimal if and only if $\sigma=\nu$.
\end{description}

\end{theorem}

\beginproof

\paragraph{1)}
The assumption (\ref {opt:theo2:assum}) guarantees that the $\nu_i(x)$ are probabilities.

\paragraph{2)}
One calculates
\be
D_\phi(\sigma||\mu)-D_\phi(\nu||\mu)
&=&
\frac 12\sum_x\left[\sigma(x)-\nu(x)\right]\,
\left[\sigma(x)+\nu(x)-2\mu(x)\right]\cr
&=&
\sum_{i=1}^n\pemp_i\frac 12\sum_{x\in O_i}\left[\sigma_i(x)-\nu_i(x)\right]\cr
& &\times
\left[\pemp_i\sigma_i(x)+\pemp_i\nu_i(x)-2p_i\mu_i(x)\right]\cr
&=&
\sum_{i=1}^n(\pemp_i)^2\frac 12\sum_{x\in O_i}\left[\sigma_i(x)-\nu_i(x)\right]^2\cr
& &
+\sum_{i=1}^n\pemp_i\sum_{x\in O_i}\left[\sigma_i(x)-\nu_i(x)\right](\pemp_i-p_i)\rho_i(x)\cr
&=&
\sum_{i=1}^n(\pemp_i)^2D_\phi(\sigma_i||\nu_i).
\nnee
In the above calculation the third line is obtained by eliminating $p_i\mu_i$
using the definition of $\nu_i$. This gives
\be
& &
\pemp_i\sigma_i(x)+\pemp_i\nu_i(x)-2p_i\mu_i(x)\cr
&=&
\pemp_i\sigma_i(x)+\pemp_i\nu_i(x)
-2\pemp_i\left[\nu_i(x)-\left(1-\frac{p_i}{\pemp_i}\right)\rho_i(x)\right]\cr
&=&
\pemp_i\left[
\sigma_i(x)-\nu_i(x)\right]
+2(\pemp_i-p_i)\rho_i(x).
\nnee
The term
\be
\sum_{i=1}^n\pemp_i\sum_{x\in O_i}\left[\sigma_i(x)-\nu_i(x)\right](\pemp_i-p_i)\rho_i(x)
\nnee
vanishes because $\rho_i(x)$ is constant on the set $O_i$
and the probability measures $\nu_i$ and $\sigma_i$ have support in $O_i$.

\paragraph{3)}
From 2) it follows that $D_\phi(\sigma||\mu)\ge D_\phi(\nu||\mu)$,
with equality when $\sigma=\nu$.

Conversely, when $D_\phi(\sigma||\mu)=D_\phi(\nu||\mu)$
then (\ref {opt:pythbreg}) implies that 
\be
\sum_{i=1}^n(\pemp_i)^2D_\phi(\sigma_i||\nu_i)&=&0.
\nnee
The empirical probabilities are strictly positive by assumption.
Hence, it follows that $D_\phi(\mu|| \sigma_i)=D_\phi(\mu||\nu_i)$ for all $i$
and hence, that $\sigma_i=\nu_i$ for all $i$. The latter implies $\sigma=\nu$.

\qed

The optimal update $\nu$ can be written as
\be
\nu=\sum_i\left[(\pemp_i-p_i)\rho_i+ p_i\mu_i\right]
=\mu+\sum_i(\pemp_i-p_i)\rho_i.
\nnee
This result is in general quite different from the update proposed by Theorem \ref{theorem1},
which is
\be
\nu&=&\sum_i\pemp_i\mu_i.
\nnee
The updates proposed by the two theorems coincide only
in the special cases that either $\pemp_i=p_i$ for all $i$
or that $\mu_i=\rho_i$ for all $i$. 
In the latter case the prior distribution $\mu=\sum_ip_i\rho_i$ is replaced by the update
$\nu=\sum_i\pemp_i\rho_i$.

The entropy of the update when event $O_i$ is observed, according to Theorem \ref{theorem1}, equals
$S(\nu_i)=S(\mu_i)$.
According to Theorem \ref{theorem2} it equals
\be
S(\nu_i)
&=&
S\left([1-\frac{p_i}{\pemp_i}]\rho_i+ \frac{p_i}{\pemp_i}\mu_i\right).
\nnee
If $p_i\le \pemp_i$ then it follows that
\be
S(\nu_i)
&\ge&
[1-\frac{p_i}{\pemp_i}]S(\rho_i)+\frac{p_i}{\pemp_i}S(\mu_i)\cr
&\ge&
S(\mu_i).
\nnee
The former inequality follows because the entropy is a concave function.
The latter follows because entropy is maximal for the uniform distribution $\rho_i$.
On the other hand, if $p_i>\pemp_i$ then one has
\be
S(\mu_i)
&=&
S\left([1-\frac{\pemp_i}{p_i}]\rho_i+ \frac{\pemp_i}{p_i}\nu_i\right)\cr
&\ge&
[1-\frac{\pemp_i}{p_i}]S(\rho_i)+\frac{\pemp_i}{p_i}S(\nu_i)\cr
&\ge&
S(\nu_i).
\nnee
In the latter case the decrease of the entropy is stronger than in the case of the update
based on the squared Hellinger distance.
In conclusion, the update relying on the quadratic Bregman divergence 
looses details of the prior distribution by making a convex combination with a
uniform distribution weighed with the probabilities of the observation.
It does this more so for the events with observed probability larger than predicted,
this is when $\pemp_i>p_i$.

Note that
Theorem \ref{theorem2}
cannot always be applied because it contains restrictions on the empirical probabilities.
In particular, if the prior probability $\mu(x)$ of some point $x$ in $X$ vanishes then
the condition (\ref {opt:theo2:assum}) requires that the empirical probability $\pemp_i$
of the partition $O_i$ to which the point $x$ belongs is larger than
or equal to the prior probability $p_i$.

\subsection{Update of conditional probabilities}
\label{sect:updcond}

The two previous theorems assume that no empirical information is available about
conditional probabilities. If such information is present
then an optimal choice should make use of it.
In one case the solution of the problem is straightforward. If the probabilities $\pemp_i$
are available together with all conditional probabilities $\pemp(B|O_i)$ and there exists
an update $\nu$ which reproduces these results then it is unique. Two cases remain:
1) The information about the conditional probabilities is incomplete; 2) the information
is internally inconsistent -- no update exists which reproduces the data.

Let us tackle the problem by considering the case that the single information that is
available besides the  probabilities $\pemp_i$ is the vector of conditional probabilities $\pemp(B|O_i)$
of a fixed event $B$ given the outcome of the measurement of the random variable $g$ as introduced in
Section \ref{sect:emp}.

The following result is independent of the choice of divergence function.

\begin{proposition}
Fix an event $B$ in $X$.
Assume that the conditional probabilities $p(B|O_i)$, $i=1,\cdots,n$, are strictly positive
and strictly less than 1.
Assume in addition that $\pemp_i\pemp(B|O_i)\le 1$ for all $i$.
Then there exists an update $\nu$ with corresponding probabilities $q(\cdot)$ such that
$q(O_i)=\pemp_i$ and $q(B|O_i)=\pemp(B|O_i)$, $i=1,\cdots,n$.
\end{proposition}

\beginproof

An obvious choice is to take $\nu$ of the form
$\nu=\sum_i\pemp_i\nu_i$ with $\nu_i$ of the form
\be
\upd\nu_i(x)=\left[a_i\Io_{B\cap O_i}(x)+b_i\Io_{\Bc\cap O_i}(x)\right]\upd\mu(x),
\nnee
with $a_i\ge  0$ and $b_i\ge 0$.
Normalization of the $\nu_i$ gives the conditions
\be
1&=&a_ip(B\cap O_i)+b_ip(\Bc\cap O_i).
\label{single:normcond}
\ee
Reproduction of the conditional probabilities gives the conditions
\be
\pemp(B|O_i)
=\frac{q(B\cap O_i)}{q(O_i)}
=a_i\frac{p(B\cap O_i)}{\pemp_i}.
\nnee
The latter gives
\be
a_i=\frac{\pemp_i}{p_i}\,\frac{\pemp(B|O_i)}{p(B|O_i)}.
\nnee
The normalization condition (\ref {single:normcond}) becomes
\be
1=\pemp_i\pemp(B|O_i)+b_ip(\Bc\cap O_i).
\nnee
It has a positive solution for $b_i$ because $\pemp_i\pemp(B|O_i)\le 1$ and $p(\Bc\cap O_i)>0$.
\qed

\subsection{The Hellinger case}
\label{sect:opt:hell}

The optimal updates can be derived easily from Theorem \ref{theorem1}.
Double the partition by introduction of the following sets
\be
O_i^+=B\cap O_i
\quad\mbox{ and }\quad
O_i^-=\Bc\cap O_i.
\nnee
They have prior probabilities $p_i^\pm=p(O_i^\pm)$.
Corresponding prior measures $\mu_i^\pm$ are defined by
\be
\upd\mu_i^\pm(x)=\frac{1}{p_i^\pm}\Io_{O_i^\pm}(x)\upd\mu(x)
\nnee

The empirical probability of the set $O_i^+$ is taken equal to $\pemp_i\pemp(B|O_i)$,
that of $O_i^-$ equals $\pemp_i[1-\pemp(B|O_i)]$.
The optimal update $\nu$ follows from Theorem \ref{theorem1} and is given by
\be
\upd\nu(x)&=&\sum_i\pemp_i\pemp(B|O_i)\,\upd\mu_i^+(x)
+\sum_i\pemp_i[1-\pemp(B|O_i)]\,\upd\mu_i^-(x).
\label{single:hell}
\ee
By construction is 
\be
q(O_i^+)=\pemp_i\pemp(B|O_i)
\quad\mbox{ and }\quad
q(O_i^-)=\pemp_i[1-\pemp(B|O_i)].
\nnee
One now verifies that $q(O_i)=\pemp_i$ and $q(B|O_i)=\pemp(B|O_i)$, which is the intended result.

\subsection{The Bregman case}

Next consider the optimization with the quadratic Bregman divergence.
Probability distributions $\rho_i^\pm$ are defined by
\be
\rho_i^\pm(x)&=&\frac{1}{|O_i^\pm|}\Io_{O_i^\pm}(x).
\nnee
Introduce the notations
\be
r_i^+&=&\frac{p_i^+}{\pemp_i\pemp(B|O_i)},\cr
r_i^-&=&\frac{p_i^-}{\pemp_i[1-\pemp(B|O_i)]},\cr
& &\cr
\nu_i^\pm(x)&=&(1-r_i^\pm)\rho_i^\pm+r_i^\pm\mu_i^\pm(x).
\nnee
Then the condition for Theorem \ref{theorem2} to hold is that $\nu_i^\pm(x)\ge 0$
for all $x,i$. The optimal probability distribution $\nu$ is given by
\be
\nu(x)&=&\sum_i \pemp_i\pemp(B|O_i)\nu^+_i(x)+\sum_i \pemp_i[1-\pemp(B|O_i)]\nu^-_i(x)\cr
&=&
\sum_i\left[\pemp_i\pemp(B|O_i)-p_i^+\right]\rho_i^+
+\sum_ip_i^+\mu_i^+\cr
& &
+\sum_i\left[\pemp_i[1-\pemp(B|O_i)]-p_i^-\right]\rho_i^-
+\sum_ip_i^-\mu_i^-\cr
&=&
\sum_i\pemp_i\pemp(B|O_i)\,\left[\rho_i^+-\rho_i^-\right]\cr
& &
-\sum_i p_i^+\rho_i^++\sum_i[\pemp_i-p^-_i]\rho_i^-\cr
& &+\mu.
\nnee

\section{Summary}

It is well known that the use of unmodified prior conditional probabilities is the optimal way
for updating a probability distribution after new data come available. The update procedure
minimizes the Hellinger distance between prior and posterior probability distributions.
For the sake of completeness a proof is given in Theorem \ref{theorem1}.

In the context of the present research the work of 
Banerjee, Guo and Wang \cite {BGW05} was considered as well. 
They prove that minimization of
the Hellinger distance can be replaced by minimization of a Bregman divergence,
without modifying the outcome.
However, their proof is restricted to random variables that are functions of
the updated random variable. It is shown in Theorem \ref{theorem2} that a more general
context yields results quite distinct from those obtained by the usual procedure.


\begin{thebibliography}{9}

\bibitem{WH82}
H. White, {\em
Maximum Likelihood Estimation of Misspecified Models,}
Econometrica {\bf  50}, 1--25 (1982).


\bibitem{BGW05}
A. Banerjee, X. Guo, H. Wang,
{\em On the Optimality of Conditional Expectation as a Bregman Predictor,}
IEEE Trans. Inf. Th. {\bf 51}, 2664--2669 (2005).

\bibitem{FSG08}
B.A. Frigyik, S. Srivastava, M.R. Gupta, {\em 
Functional Bregman Divergences and Bayesian Estimation of Distributions,}
IEEE Trans. Inf. Th. {\bf 54}, 5130--5139 (2008).

\bibitem{AN00}
S. Amari, H. Nagaoka, {\em
Methods of Information Geometry}
(Oxford University Press, 2000)
(Originally published in Japanese by Iwanami Shoten, Tokyo, Japan, 1993) 

\bibitem{AS16} S. Amari, {\em
Information Geometry and its Applications}
(Springer, 2016)

\bibitem{AJLS18}
N. Ay, J. Jost, H. V\^an L\^e, L. Schwachh\"ofer, {\em
Information Geometry} (Springer, 2017).

\bibitem{JE57}
E. Jaynes, {\em
Information theory and statistical mechanics,}
Phys. Rev. {\bf 106}, 620--630 (1957).

\end{thebibliography}
\end{document}